\documentclass[11pt] {article}      

\usepackage[top=.9in, bottom=.9in, left=1.6in, right=1.6in]{geometry}

\usepackage{amsmath,amssymb,eucal}
\usepackage[T1]{fontenc}  
\usepackage{graphicx}   
\usepackage{enumerate}  
\usepackage{url}               
\usepackage{float}  
\usepackage{tikz}
\usetikzlibrary{matrix, arrows}
\usetikzlibrary{calc}

\usepackage{newtxtext,newtxmath}

\def\g{g}                               
\newtheorem{theorem}{Theorem}[section]            
\newtheorem{proposition}[theorem]{Proposition}  
\newtheorem{corollary}[theorem]{Corollary}	      %
\numberwithin{equation}{section}
\numberwithin{figure}{section}
\newcommand{\gaction}[2]{\genfrac{}{}{0.5pt}{}{#1}{#2}%
                        \!\lower2pt\hbox{\rotatebox[origin=c]{-90}{{$\looparrowright$}}}}

\newcommand{\dotfraction}[2]{\genfrac{}{}{0.5pt}{}{#1}{#2}%
                        \!\lower.5pt\hbox{{$\circ$}}}

\def\PSL{\rm PSL}
\def\SU{\rm SU}
\def\PSU{\rm PSU}

\begin{document}
\title{On the Dedekind tessellation}

\author{Jerzy Kocik \\ 
\small Department of Mathematics\\
\small Southern Illinois University\\
\small jkocik{@}siu.edu}

\date{}

\maketitle

\begin{abstract}
The Dedekind tessellation --  the regular tessellation of the upper half-plane 
by the M{\"o}bius action of the modular group -- is usually 
viewed as a system of ideal triangles.
We change the focus from triangles to circles and give their complete 
algebraic characterization with the help of a representation 
of the modular group acting by Lorentz  transformations on Minkowski space. 
This interesting example of the interplay of geometry, group theory and number theory 
leads also to convenient algorithms for computer drawing of the Dedekind tessellation.
\\[7pt]
{\small 
}
\end{abstract}

\section{Introduction}

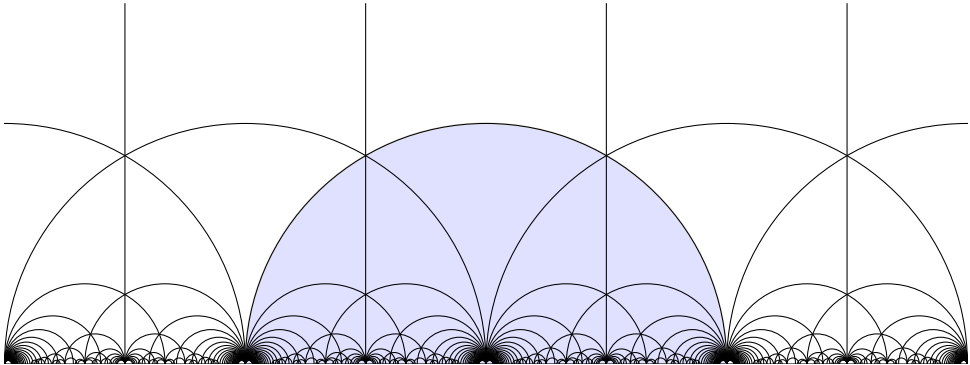
\begin{figure}[H]
\centering
\begin{tikzpicture}[scale=3.2]
\clip (-2,0) rectangle (2,1.5);
\draw [fill=blue, opacity=.125, thick] (0,0) circle (1);
\draw (2.5,0)--(-2.5,0) ;   
\foreach \n in {0, 1, 2, 3}
\draw (\n + 0.5, -1.5)--(\n + 0.5, 1.5) 
          (-\n - 0.5, -1.5)--(-\n - 0.5, 1.5); 
\foreach \n in {-2, -1, 0, 1, 2}
\draw (\n , 0) circle (1/1);
\foreach \n in {-3,-2, -1, 0, 1,2, 3}
\foreach \a/\b   in {  
1/3,  1/5,  1/7,  3/8,  1/9,  1/11,  1/13,  1/15,  4/15,  7/16,  1/17,  1/19,  1/21,  8/21,  1/23,  5/24,  11/24,
1/25,  1/27,  1/29,  1/31,  15/32,  1/33,  1/36,  6/36,  1/37,  1/39,  14/39,  11/40,  19/40,  1/41,  1/43,
1/45,  19/45,  1/47,  7/48,  23/48,  1/49,  1/51,  16/51,  1/53,  1/55,  21/55,  13/56,  27/56,  1/57,  20/57,
1/59, 1/61, 1/63, 8/63, 31/64, 1/65, 14/65, 1/67, 1/69, 22/69 
}
\draw (\n + \a/\b, 0) circle (1/\b)
           (\n - \a/\b, 0) circle (1/\b)
;
\end{tikzpicture}
\caption{The Dedekind tessellation}
\label{fig:wiki}
\end{figure}

Figure (\ref{fig:wiki}) shows the well-known {\it regular tessellation} 
of the Poincar\'e half-plane $\mathbf H$, 
the upper half of the plane of complex numbers. 
It was first discussed in the 1877 paper by Richard Dedekind \cite{Ded}, 
which justifies the name {\it Dedekind tessellation},
but the first illustration appeared in Felix Klein's paper \cite{Kn} 
(see \cite{LeB,Sti} for these credits).  
Usually, the tessellation is understood in the context of the action of 
the modular group $\PSL^\pm(2,\mathbb Z)$  on $\mathbf H$ 
as a model of non-Euclidean, hyperbolic geometry. 
It is viewed as a set of "ideal triangles'' (triangles with arc sides) 
that results as an orbit of any of them via the action of the modular group.

~

\noindent
Recall that the group $\PSL(2,\mathbb Z)$ acts on the extended complex plane 
$\dot {\mathbb C} = {\mathbb C} \cup \{\infty\}$ (Riemann sphere)
via M\"obius transformations defined by
\begin{equation}
\label{eq:moebius}
	     z \ \mapsto \ \left[ \begin{array}{cc}
                                             a & b \\
                                             c & d \end{array} \right]
                              \cdot z
                            \ \equiv \ \frac{az+b}{cz+d}\,,
\qquad
                           ad - bc = 1.  
\end{equation}
with $a,b,c,d\in\mathbb Z$.
The group preserves the Poincar\'e half-plane 
and in particular the extended real line $\dot{\mathbb R}=\mathbb R\cup\{\infty\}$.  
The Dedekind tessellation coincides with  
the orbit $\PSL(2,\,\mathbb Z) \cdot\Delta$,  where $\Delta$ is any of the ideal triangles,
but it is customary to choose 
the triangle ($|z|^2 \geq 1$,   $|\hbox{Re}(z)| \leq \frac{1}{2}$),
called the {\it fundamental domain}, 
shown shaded in Figure \ref{fig:symbols}).       
\\
\\
One may ask what would be an efficient, practical,  way 
to draw this pattern. (Producing the orbit of the modular group is hardly convenient,
since the algorithm walks over circles in an inconsistent way.)
\\
\\
{\bf Problem:} Suppose we want to view the tessellation as a system of {\bf circles} 
rather than ideal triangles.  Is there any simple characterization of 
this system of ``Dedekind circles'' --- denoted $\mathcal D$ --- in terms of the list of the centers 
(the values of $x\in\mathbb R$) and the associated radii?  
Here it is:
\\
\\
{\bf The main result:}  Algebraic characterization of the system of 
Dedekind circles: 

\begin{enumerate}
\item\vskip-5pt
The centers of the circles are rational and the curvatures of all circles are integral.  
If the center is $k/n$ then the radius is $1/n$. 
\item\vskip-3pt
For a fraction $k/n$ to be the center of a Dedekind circle, 
the following necessary and sufficient conditions must be satisfied:
\begin{itemize}
\item[$\circ$]\vskip-2pt
        		either $n$ is odd and $n\,|\,(k^2-1)$ 
\item[$\circ$]\vskip-2pt
         	or $n$ is a multiple of $8$ and $(k^2-1)/n$ is an odd integer
\item[$\circ$] \vskip-2pt
       		or $n = 0$, in which case the circle is a vertical line at $x = k/2$,  $k$ odd.
\end{itemize}
\item\vskip-3pt
The above rules exhaust the system of Dedekind circles.
\end{enumerate}

\noindent
A simple algorithm to build the set of circles of arbitrary size emerges: 
For every curvature $n$ (denominator) there is a finite set $K(n)$ of integers $0 \leq k < n$ 
such that $k/n$  is a center of a Dedekind circle (or a line if $n=2$).  
These data correspond to circles in the interval $[ 0,1 )$.
The rest of the pattern is reconstructed by the translational symmetry, 
since $z\to z+1$ is an element of $\PSL(2,\,\mathbb Z)$.
\\
\\
Table \ref{fig:table} shows the data for curvatures up to $n=107$.  
For instance the entry  ``$n=72$, $K(n)= 19, 35, 37, 53$" means that the 
interval $[0, 1]$ contains four circles of radius $r = 1/72$ with centers 
at  $x = 19/72$,  $x = 35/72$, $x = 37/72$, and $x = 53/72$.  
In Figure \ref{fig:mars}, these fractions $k/n$ are drawn as points $(k,n)$ in 
$\mathbb N\times\mathbb N$.

{
\begin{figure}[H]

\centering
\includegraphics[scale=.77]{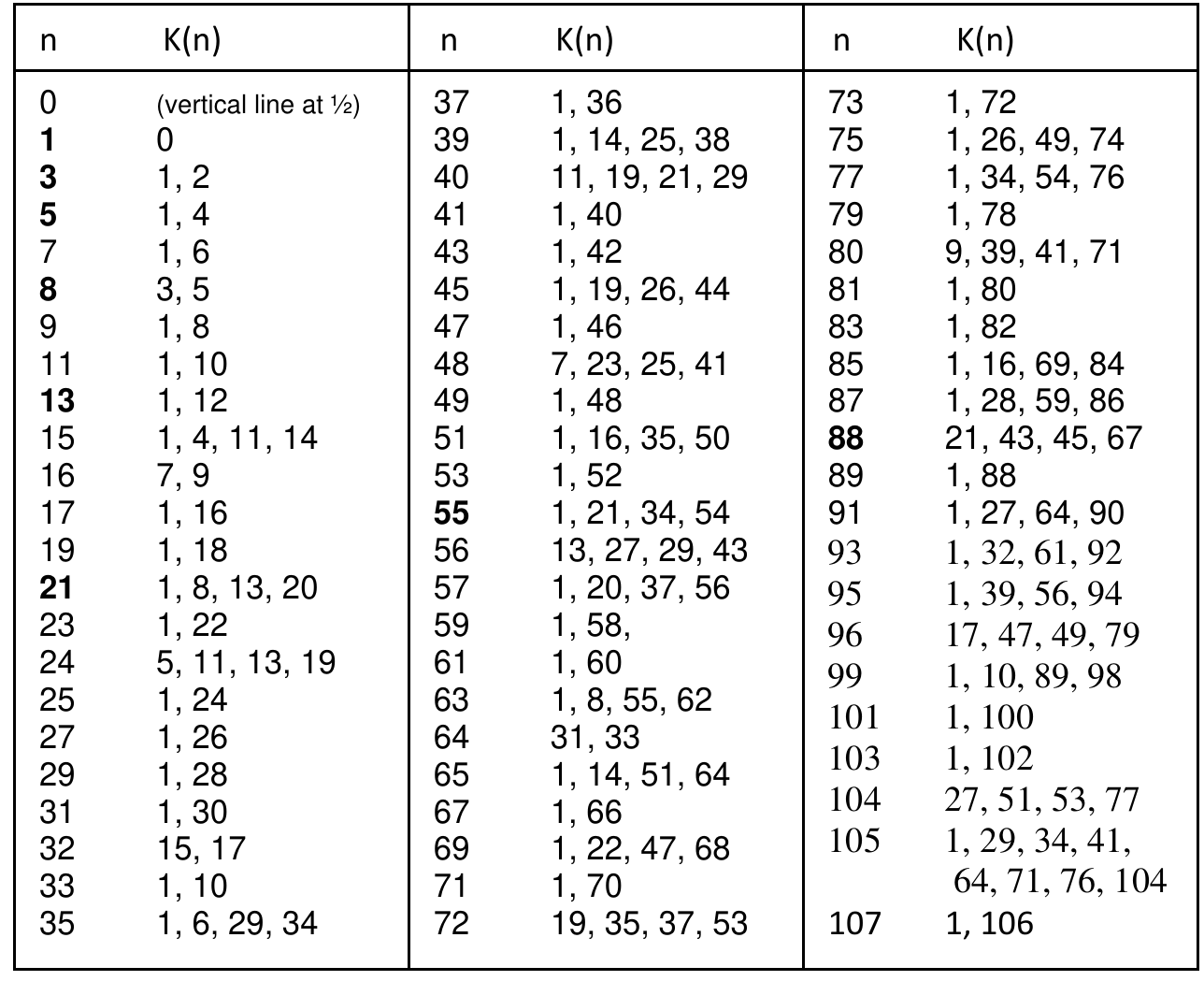}
\caption{\small Table of admissible circles with curvatures $n<108$}
\label{fig:table}
\end{figure}

\begin{figure}[H]
\centering
\includegraphics[scale=.85]{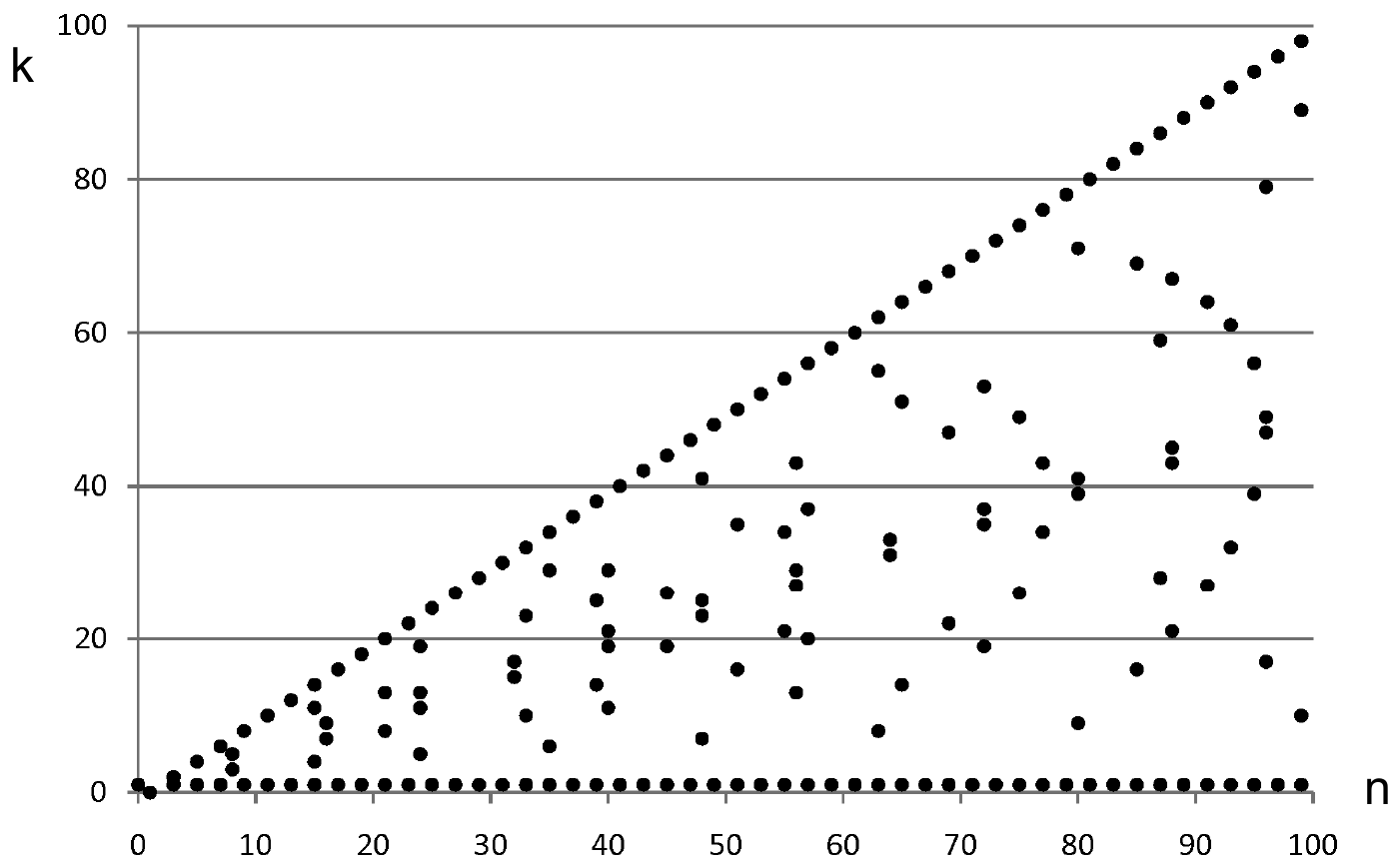}
\caption{\small ``Canals on Mars" --- which of the visible patterns are accidental 
and which have number-theoretic significance?  
Fractions $k/n$ representing the positions of the centers of the Dedekind circles
are presented as points $(n,k)$.  Only $k<n$ are shown.}
\label{fig:mars}
\end{figure}
}

A geometric object that comes to mind is the Apollonian Window \cite{jk2},  
an Apollonian gasket made of circles whose curvatures are also integral \cite{LMW, jk-Descartes}.  
The intriguing algebraic-geometric interactions in the Apollonian disk packing 
occur also in the case of the Dedekind system, and we shall use similar 
tools to investigate them and to prove the main result.

\section{Dedekind circles in Minkowski space}

Circles (and lines as their special case) may be mapped to 
the Minkowski space $\mathbf M \cong {\mathbb R}^{3,1}$, 
where their image lies on the unit space-like hyperboloid
\cite{jk-Descartes}. 
In particular, a circle centered at $(x, y)$ with radius $r$ is mapped to the vector 
\begin{equation}
C  =   \left[ \begin{array}{c}
                      \dot  x \\
                      \dot  y \\      
                       \beta \\
                       \gamma  \end{array}\right] 
\end{equation} 
where $\beta = 1/r$ is the curvature,  and $\dot x= x/r$  and $\dot y= y/r$ are ``reduced coordinates".  
The last entry, ``co-curvature" $\gamma$, has the geometric meaning of 
curvature of the image of the circle under inversion in the unit circle at the origin.    
In the case of a line, the curvature vanishes, $\beta=0$, and the correspondence is
\begin{equation}
C  =   \left[ \begin{array}{c}
                      \dot  x \\
                      \dot  y \\      
                       0 \\
                       \gamma  \end{array}\right] 
\qquad\Rightarrow\qquad
\hbox{Line} \ \ 2\dot x x + 2 \dot y y = \gamma
\end{equation} 
In general, the co-curvature may be easily determined 
from the fact that the vector is normalized, $|C|^2 = -1$, 
where the quadratic form in $\mathbf M$ is defined by
\begin{equation}
\label{eq:normbig}
\langle C,\, C\rangle = |C|^2 = -\dot x^2 - \dot y^2 +\beta\gamma
\end{equation} 
(non-diagonalized basis), i.e.  $\gamma = (\dot x^2 +\dot y^2 -1)/\beta$.   
We shall also use a more suggestive representation of the circle vectors 
that we call circle symbols \cite{jk-Descartes}, namely
\begin{equation}
 \left[ \begin{array}{c}
                      \dot  x \\
                      \dot  y \\      
                       \beta \\
                       \gamma  \end{array}\right] 
        \  = \  \frac{\dot x,\, \dot y}{\beta,\, \gamma}
        \ = \ \frac{\dot x,\, \dot y}{\beta}
\end{equation} 
where the co-curvature may optionally be dropped in case of a proper circle (not line).  
The true position of the circle may be read off directly from the symbol:
\begin{equation}
\frac{\dot x,\, \dot y}{\beta}
\qquad\Rightarrow\qquad
x=\frac{\dot x}{\beta}, \ \
y=\frac{\dot y}{\beta}, \ \
r=\frac{1}{\beta} \ 
\end{equation} 
The last comment concerns the sign convention: it is better to regard  
$C$ as a disk rather than a mere circle.  Then vector $(-C)$ 
describes the external region bounded by the same circle.  
It will have a negative curvature.  For example:                      
\begin{equation}
\label{eq:example}
\begin{array}{cl}
\frac{3,\,4}{6}           & \ \Rightarrow \ \hbox{\small disk centered at (1/2, 2/3) with radius 1/6} \\ [4pt]   
\frac{30,\,40}{60}     &\ \Rightarrow \ \hbox{\small disk centered at (1/2, 2/3) with radius 1/60} \\ [4pt] 
\frac{-3,\,-4}{-6}       &\ \Rightarrow \  \hbox{\small complement of the first disk (center at (1/2, 2/3), radius -1/6)}    
\end{array}
\end{equation}
Similarity transformations of the Euclidean space induce transformations of circles, 
and these induce Lorentz transformations in the Minkowski space $\mathbf M$. 
This defines a group homomorphism from the conformal group on a plane to the Lorentz ${\rm O}(3,1)$. 
Since they do not seem to appear in the literature, a number of transformation 
correspondences are given below explicitly:

\begin{proposition}
\label{prop:transf}
The basic transformations of the Euclidean plane, 
translations, dilations, inversions, reflections and rotations 
and their compositions induce the Lorentz transformations in the Minkowski space. 
In particular:
$$
\begin{array}{lcrl}
            T_v = \hbox{translation by}\  v = (a,b)                                                  
&&   \mathbf T_v 
&=  \left[ \begin{array}{cc|cc}
                          1&0&a&0 \\
                          0&1&b&0 \\
\hline
                          0&0&1&0 \\
                         2a&2b&a^2+b^2&1 \\
 \end{array}\right] 
\\[25pt]
          N = \hbox{inversion in the unit central circle}  
&& 
            \mathbf N 
&=  \left[ \begin{array}{cc|cc}
                          1&0&0&0 \\
                          0&1&0&0 \\ \hline
                          0&0&0&1 \\
                         0&0&1&0 \\
 \end{array}\right] 
\\[25pt]
           D_s = \hbox{dilation by}\ s \ \hbox{with respect to the origin}            
&& 
            \mathbf D_s 
&=  \left[ \begin{array}{cc|cc}
                          1&0&a&0 \\
                          0&1&0&0 \\ \hline
                          0&0&1/s&0 \\
                         0&0&0&s \\
 \end{array}\right] 
\end{array}
$$

\end{proposition}
\noindent
{\bf Proof:}  The proof is calculational: simply convert the vectors to 
the Euclidean description, apply a transformation, and convert 
back to the Minkowski space.  Denote the circle centered at $(x,y)$ with 
radius $r$ by $C(x, y; r)$.  Then, for instance for the translation, we have:
$$
\begin{array}{lll}
 \left[ \begin{array}{c}
                      \dot  x \\
                      \dot  y \\      
                       \beta \\
                       \gamma  \end{array}\right] 
\ \to \ 
 C(\frac{\dot x}{\beta},\, \frac{\dot y}{\beta}; \, \frac{1}{\beta})
\ \xrightarrow{~~~T~~} \ 
&
 C(\frac{\dot x}{\beta}+a,\, \frac{\dot y}{\beta}+b,\, ; \frac{1}{\beta})
\\[-30pt]
&=  C(\frac{\dot x+a}{\beta},\, \frac{\dot y+b}{\beta};\, \frac{1}{\beta})
&\ \to\  \left[ \begin{array}{c}
                      \dot  x +a \\
                      \dot  y + b\\      
                       \beta \\
                       \gamma '  \end{array}\right] 
\end{array}
$$
where the new $\gamma'$ may be easily calculated:
$$
\begin{array}{ll}
\displaystyle \frac{(\dot x + a\beta)^2 + (\dot y + b\beta)^2  - 1}{\beta}
 & = \displaystyle \frac{(\dot x^2 + \dot y^2  - 1) + 2\dot x a\beta  + 2\dot y b\beta +(a^2+b^2)\beta^2}{\beta}   \\[10pt]
 & = \gamma + 2\dot x a\beta  + 2\dot y b\beta +(a^2+b^2)\beta
\end{array}
$$
 The metric of the Minkowski space is 
$$
\g  =   \left[ \begin{array}{cc|cc}
                         -1&0&a&0 \\
                          0&-1&0&0 \\ \hline
                          0&0&0&1/2 \\
                         0&0&1/2&0
 \end{array}\right] 
$$
It is easy to verify that any transformation $A$ from the above list 
satisfies $A^T \g A = \g$, that is it preserves the metric and therefore the norm.  
$\square$

~

\begin{figure}[h]
\centering
\includegraphics[scale=.8]{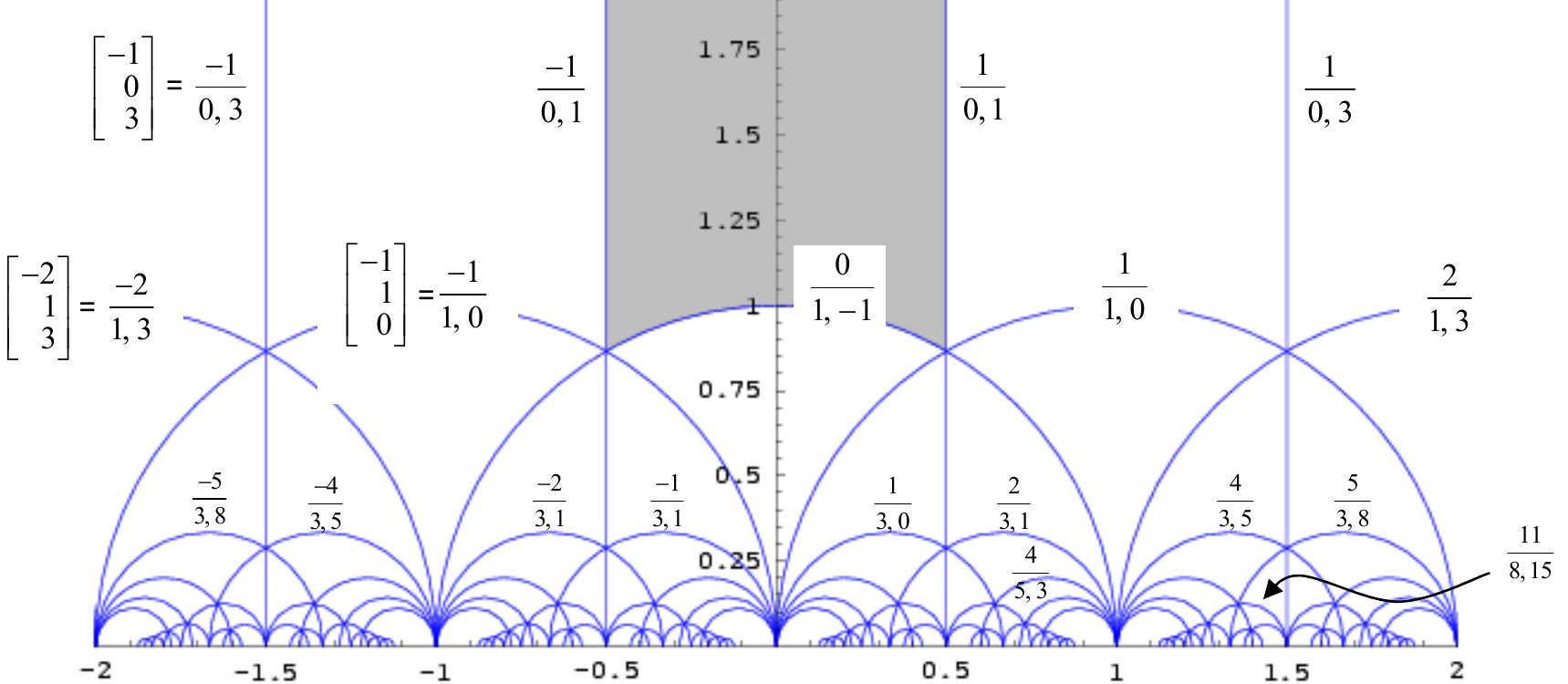}
\caption{\small Dedekind circles labeled with their representations in the Minkowski space}
\label{fig:symbols}
\end{figure}

\noindent
Now let us turn to the Dedekind system of circles.  
Since their centers are located on the real axis, $y = 0$, 
we shall omit it in the symbol and vector description:  
\begin{equation}
 C \ = \ \left[ \begin{array}{c}
                      \dot  x \\
                         0  \\      
                       \beta \\
                       \gamma  \end{array}\right] 
    \ \mapsto \  \left[ \begin{array}{c}
                      \dot  x \\
                      \beta \\
                       \gamma  \end{array}\right] 
    \ \equiv \  \frac{\dot x}{\beta,\, \gamma}
    \ \mapsto  \ \frac{\dot x}{\beta} \ \mapsto \ x
\end{equation} 
where $\beta = 1/r$, $\dot x = x/r$.   
In effect, the Dedekind circles form a 3-dimensional Minkowski 
subspace $\mathbb R^{2,1} \subset \mathbb R^{3,1} = M$.  
In the special case of vanishing curvature, $\beta=0$, the vector represents the line
\begin{equation}
\left[ \begin{array}{c}
                      \dot  x \\
                      0 \\
                       \gamma  \end{array}\right] 
 =  \frac{\dot x}{0,\, \gamma}
\quad\Rightarrow\quad
\hbox{vertical line at } \  
x = \frac{\gamma}{2\dot x}
\end{equation} 
Figure \ref{fig:symbols} shows these symbols in the Dedekind tessellation.  
The whole idea is summarized in the following diagram:
\def\smalll{\scriptsize}
\begin{equation}
\begin{tikzpicture}[baseline=-0.8ex]
    \matrix (m) [ matrix of math nodes,
                         row sep=2em,
                         column sep=4em,
                         text height=3.8ex, text depth=3ex] 
   {
   \displaystyle\gaction{\ \ \PSL(2,\mathbb C)\ }{\dot{\mathbb C}}  \quad   
                    & \quad \displaystyle\gaction{\ \ \PSL^\pm(2,\mathbb C )\ }{\mathcal C}  \quad  
                          & \quad \displaystyle\gaction{\ \ {\rm O}_{3,1}(\mathbb R) \  }{  \mathbb R^{3,1}}   \\
    \displaystyle\gaction{\ \ \PSL(2,\mathbb Z)\ }{\mathbf H} \quad
                    & \quad \displaystyle\gaction{\ \ \PSL^\pm(2,\mathbb Z )\ }{\mathcal D} \quad
                          & \quad \displaystyle\gaction{\ \ {\rm O}_{2,1}(\mathbb Z) \  }{\mathbf M \cong \mathbb R^{2,1}}       \\
     };
    \path[->]
        (m-1-1) edge node[above] {$\hbox{ \sf\smalll  induce}$} (m-1-2)
        (m-1-2) edge node[above] {$\hbox{ \sf\smalll  homo}$} (m-1-3)
        (m-2-1) edge node[above] {$\hbox{ \sf\smalll  induce }$} (m-2-2)
        (m-2-2) edge node[above] {$\hbox{ \sf\smalll  homo }$} (m-2-3)
        (m-2-1) edge node[right]  {$\imath$} (m-1-1)
        (m-2-2) edge node[right] {$\imath$} (m-1-2)
        (m-2-3) edge node[right] {$\imath$} (m-1-3);
\end{tikzpicture}   
\end{equation}

The fraction-like symbol $\gaction{\quad}{}$ denotes a group (numerator) and the set  on which it acts (denominator).
The top row reads: the modular group $\PSL(2,\mathbb C)$ may be extended to the full M\"obius group;  
its point-wise action induces its action on the circles $\mathcal C$ (which include 
lines as special cases).  This in turn may be viewed formally as an action of the Lorentz group (orthogonal group)
on the Minkowski space, from which only the unit space-like hyperboloid is the image of $\mathcal C$.
The bottom row is a reduction of these to (1) circles with centers on the real line 
and to (2) the ring $\mathbb Z$  (symbols $\imath$ at the vertical arrows denote the inclusions). 
\\

\newpage
\section{Proof of the main result} 

\begin{figure}[h]
\centering
\includegraphics[scale=1]{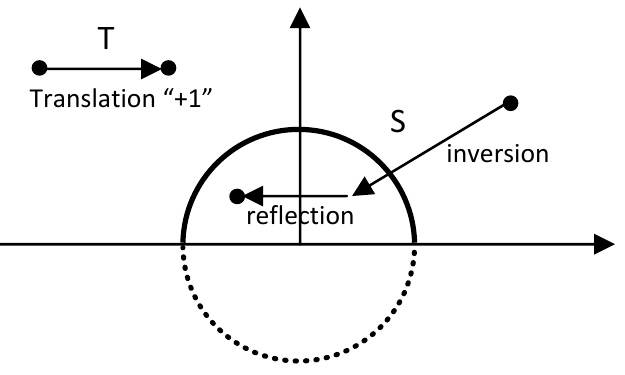}
\caption{\small The geometry of generators $T$ and $S$ of the modular group}
\label{fig:action}
\end{figure}

As is well known, 
the modular group can be generated by two elements, 
for instance $\PSL(2,\,\mathbb Z) = \hbox{gen}\{T,S\}$, where
\begin{equation}
\label{eq:TS}
\begin{array}{rccc}
         \hbox{ translation}     &T: z\to z+1.   &\qquad\hbox{As the matrix,}  &
        T=\left[ \begin{array}{cc}
                      1&1 \\
                      0 & 1  \end{array}\right]                          \\[12pt]
           \hbox{quasi-inversion}    &S: z\to -\frac{1}{z}.   &\qquad\hbox{As the matrix,}  &
         S=\left[ \begin{array}{cr}
                      0&-1 \\
                      1 & 0  \end{array}\right]                          
\end{array}
\end{equation} 
The first element is a unit translation along the $x$-axis. 
The second amounts to an inversion through the unit 
circle at the origin followed by reflection in the y-axis (Figure \ref{fig:action}).  
Although not exactly a pure inversion, it is often called so for 
terminological simplicity; we shall prefer `quasi-inversion' for clarity.

The generators satisfy the following relations:
\begin{equation}
\label{eq:identities}
                T^2 = \hbox{id}   \quad \hbox{and}\quad     (ST)^3 = \hbox{id}
\end{equation}
The group $\PSL(2,\mathbb Z)$ is covered by the algebra $\{\,T, S\,\}^*$ of words over the alphabet $\{\,T,\, S\,\}$.  
For instance the inverse translation may be seen as $T^{-1} = STSTS$, due to (\ref{eq:identities}).  

\begin{proposition} 
\label{prop:TS}
(Geometrization of the problem)  The system $\mathcal D$ of circles is the orbit of the central unit 
circle $C_0$  via action of the group
\begin{equation}
\begin{array}{ccl}
 			\mathcal D   &=&   \PSL(2,\, \mathbb Z) \cdot C_0 \\
 				      &=&  \{T, \,S\}^* \cdot C_0 
\end{array}
\end{equation}
\end{proposition}

~\\
{\bf Proof:} Indeed, the boundary of the fundamental region consists of 
pieces of the images of $C_0$, namely $C_0$, $STC_0$ and $TSTC_0$. 
$\square$

~\\
Let $N$ denote inversion in the unit central circle and let $F$ 
denote the reflection through the $y$-axis, $z \to -\bar z$.  Clearly $S = FN = NF$.  
Neither $F$ nor $N$ are in $\PSL(2,\, \mathbb Z)$ as they are odd elements of the extended group 
of all conformal transformations, $\PSL^\pm(2,\,\mathbb Z)$.  
But since we are interested  in the geometric aspect of $\mathcal D$, we shall use them as follows:

\newpage

\begin{proposition}  %
\label{prop:TI}
The Dedekind system of circles $\mathcal D$ coincides with the orbit
$$
\mathcal D = \{\,T,\, T^{-1},\, S \,\}^* \cdot C_0
$$
\end{proposition}  

\noindent
{\bf Proof:} 
Consider $C\in\mathcal D$. Starting with Proposition \ref{prop:TS}, we have the following development:
$$
\begin{array}{ccllc}
 		C   &=&  w\cdot C_0      &\quad    \hbox{for some} \  w\in \{\,S,\, T\,\}^*        &\qquad (i) \\
 		     &=&  w' \cdot C_0    &\quad    \hbox{for some} \  w'\in \{\,N,\,F,\, T\,\}^*   &\qquad (ii) \\     
 		     &=& w" \cdot C_0    &\quad    \hbox{for some} \  w"\in \{\,N,\,T,\,T^{-1}\,\}^*   &\qquad (iii) 
\end{array}
$$
where $(ii)$ results by replacing every $S$ in $w$ by $NF$,  
and $(iii)$ follows by pushing every $F$ in $w'$ to the right using the following 
commutation relations
\begin{equation}
(a) \   FT = T^{-1}F   \qquad  (b) \   NF = FN
\end{equation}
Now we get $C = w" F^n \cdot C_0$ for some $n\in \mathbb N$.  But $F$ leaves $C_0$ invariant, 
hence the result. $\square$
\\

Recall that the Lorentz transformations of the Minkowski space $\mathbf   M$
corresponding to elements of the modular group $\PSL(2,\mathbb Z)$ 
as acting on the space of circles are as follows (vectors written as circle symbols):
\begin{equation}
\label{eq:transfbold}
\begin{array}{ll}
 \mathbf N  : \displaystyle\frac{\dot x}{\beta, \, \gamma} \ \to \  \frac{\dot x}{\gamma, \, \beta}   
                    &   \quad\quad      \mathbf N^{-1} = \mathbf N   \\ [10pt]
    	\displaystyle \mathbf T :   \frac{\dot x}{\beta, \, \gamma} \ 
                                                                 \to \  \frac{\dot x+\beta}{\beta, \, \beta+2\dot x+ \gamma}   
           &    \displaystyle\quad\quad     \mathbf T^{-1}:  \frac{\dot x}{\beta, \, \gamma} \ 
                                                                 \to \  \frac{\dot x-\beta}{\beta, \, \beta-2\dot x+ \gamma}   
\end{array}
\end{equation}
Recall also that the normalization, expressed as
\begin{equation}
\label{eq:norm}
\dot x^2 - 1 = \beta\gamma
\end{equation}
is preserved by the Lorentz transformations including (\ref{eq:transfbold}).  
Now we may restate Proposition \ref{prop:TI}:

\begin{corollary} 
\label{prop:ITD}
The Dedekind system of circles coincides with the words over the 
alphabet $\{\,N,\,T,\,T^{-1}\}$ applied to the single vector representing the central circle
$$
    	\mathcal D \ = \ \left\{\,\mathbf N,\,\mathbf T,\,\mathbf T^{-1}\right\}^* \cdot \left(\frac{0}{1,-1}\right)
$$
\end{corollary}

\begin{proposition}  
\label{prop:gcd}
Any symbol in $\mathcal D$ satisfies $\gcd(\dot x,\beta) =\gcd(\dot x,\gamma) = 1$,  
unless $\beta$ or $\gamma$ equal zero.
\end{proposition} 

\noindent
{\bf Proof:}  Suppose not. Say $\dot x = kp$ and $\beta = kq$ for some integer $k>0$.  
Then (\ref{eq:norm}) gives $(kp)^2 -1 = kq\gamma$, 
or  $kq\gamma = (kp -1)(kp+1)$, 
which is impossible:  $k\,|\,(kp-1)$ or $k\,|\,(kp+1)$ only if $k = 1$. 
The same goes for $\gamma$. 
$\square$

\begin{proposition}  
All circles of the Dedekind tessellation have 
centers at rational numbers on the real axis, 
and if the center is at $k/n$, than the radius is $1/n$. 
\end{proposition} 

\noindent
{\bf Proof:}   
Matrices $\mathbf T$ and $\mathbf N$ preserve the integrality of $\dot x$ and $\beta$, 
hence, via Corollary~\ref{prop:ITD}, $x = \dot x /\beta\in\mathbb Q$  and $r = 1/\beta$.  
(refer to the example in \ref{eq:example} as to the importance of Proposition \ref{prop:gcd}.
$\square$
\\

Thus the problem reduces now to identifying the fractions representing the centers of the circles of $\mathcal D$.  
Any integer curvature $n$ has the associated set of numbers  $K(n) = {k_1,k_2, \ldots}$  
that are permitted to serve as the numerator in $k/n$.  
We assume that $k<n$, since the circles with centers in $[0,1)$ are sufficient to recover the whole tessellation.

We will need yet another observation. Define the parity of a set $A\subset \mathbb Z$ as the 
pair $(m,n)$ where $m$ denotes the cardinality of the odd members in $A$ and $n$ the even ones.  
For example, the parity of $\{0, 1, 3, 7, 8\}$ is $(3,2)$.

\begin{proposition}  
\label{prop:parity}
In the symbol of any Dedekind circle, exactly two of the three numbers are odd; 
in other words, the parity of $\{\dot x,\beta,\gamma\}$ is $(2,1)$.
\end{proposition} 

\noindent
{\bf Proof:} Note that $\mathbf T$ and $\mathbf N$ of (\ref{eq:transfbold}) change 
the parity of the symbol as follows:
\begin{gather*}
{\scriptstyle \mathbf N} \  
\hbox{\Large\rotatebox[origin=c]{-90}{$\circlearrowright$}}\ \
\frac{\hbox{\small even}}{\hbox{\small odd},\hbox{\small odd}}
\ \xleftarrow{\ \ \mathbf T\!}\!\!\to \
\frac{\hbox{\small odd}}{\hbox{\small odd},\hbox{\small even}}
\ \xleftarrow{\ \ \mathbf N}\!\!\to \
\frac{\hbox{\small odd}}{\hbox{\small even},\hbox{\small odd}}
\ \ \hbox{\Large\rotatebox[origin=c]{90}{$\circlearrowleft$}}\ {\scriptstyle \mathbf T}
\\[6pt]
{\scriptstyle \mathbf N} \  
\hbox{\Large\rotatebox[origin=c]{-90}{$\circlearrowright$}}\ \
\frac{\hbox{\small odd}}{\hbox{\small odd},\hbox{\small odd}}
\ \xleftarrow{\ \ \mathbf T\!}\!\!\to\
\frac{\hbox{\small even}}{\hbox{\small odd},\hbox{\small even}}
\ \xleftarrow{\ \ \mathbf N}\!\!\to\
\frac{\hbox{\small even}}{\hbox{\small even},\hbox{\small odd}}
\ \ \hbox{\Large\rotatebox[origin=c]{90}{$\circlearrowleft$}}\ {\scriptstyle \mathbf T}
\\[6pt]
{\scriptstyle \mathbf N} \  
\hbox{\Large\rotatebox[origin=c]{-90}{$\circlearrowright$}}\ \
\frac{\hbox{\small odd}}{\hbox{\small even},\hbox{\small even}}
\ \ \hbox{\Large\rotatebox[origin=c]{90}{$\circlearrowleft$}}\ {\scriptstyle \mathbf T}
\qquad\qquad\qquad
{\scriptstyle \mathbf N} \  
\hbox{\Large\rotatebox[origin=c]{-90}{$\circlearrowright$}}\ \
\frac{\hbox{\small even}}{\hbox{\small even},\hbox{\small even}}
\ \ \hbox{\Large\rotatebox[origin=c]{90}{$\circlearrowleft$}}\ {\scriptstyle \mathbf T}
\end{gather*}
\\
(four closed orbits).  In particular, both $\mathbf T$ and $\mathbf N$ preserve parity $(2,1)$, 
as shown in the first row.  Since every symbol in $\mathcal D$ is an image of $0/(1,-1)$, 
the central circle, which is of parity $(2,1)$, the conclusion follows.  
$\square$

\begin{theorem}  
\label{theorem:main}
The Dedekind system $\mathcal D$ is composed of the following elements:
\begin{itemize}
\item[$(A)$]\vskip-5pt
Circles of curvature $n$ at every rational number $x=k/n$ satisfying one of the following conditions: \\[-17pt]
\begin{itemize}
          \item[$(i)$ ~]  \vskip-3pt     $n$ is odd and $n\,|\,(k^2-1)$

          \item[$(ii)$ ~]  \vskip-3pt     $n$ is a multiple of eight 
                                                    and $(k^2-1)/n$ is an odd integer
\end{itemize}
\item[$(B)$]\vskip-5pt
         	 Lines at every $k/2$,      ($k$ odd)
\end{itemize}
\end{theorem}

\noindent
{\bf Proof:} ($\Rightarrow$)
First, we show that the elements of $\mathcal D$ must satisfy the above conditions. 
Indeed, set $\dot x/(\beta,\gamma)\in\mathcal D$.   
All these entries are integer due to Corollary \ref{prop:ITD} and integrality 
of transformations (\ref{eq:transfbold}). 

Suppose $\beta\not= 0$.  
Now $\dot x$  is either odd or even.  
If $\dot x$  is even, then the first condition ({\it i}\/) of $(A)$ 
follows by Proposition \ref{prop:parity}  and normalization (\ref{eq:norm}).  
If it is odd ($\dot x = 2p+1$ for some $p$), then the normalization (\ref{eq:norm}) implies:
$$
                   (2p+1)^2 - 1 = \beta\gamma \,,
$$
that is:
$$
                 4p(p+1) = \beta\gamma  \,.
$$
But $p(p+1)$ is divisible by $2$ and is equal to $2\Delta$ for some triangular integer $\Delta$.
Therefore:
$$
                    \beta\gamma = 8\Delta\,.
$$
Hence $8\,|\,\beta\gamma$.  But 8 must divide only one of the factors, $\beta$ or $\gamma$, 
because only one can be even due to Proposition \ref{prop:parity}. 
In conclusion, if $2\,|\,\beta$  then $8\,|\,\beta$, and $\dot x$ must be odd.  
This gives the second case of (A).  

In the case of $\beta=0$, the symbol must be of the form 
$$
           \frac{1}{0,\gamma}
$$
due to the norm condition (\ref{eq:norm}),  where $\gamma$ is odd due to Proposition \ref{prop:parity}.  
This describes the vertical line at $x = \gamma/2$, noted as case (B)  
in the theorem.

~\\
($\Leftarrow$) \   Now we shall see that any symbol that satisfies the conditions of the theorem is 
an element of $\mathcal D$.  The argument is much like Euclid's algorithm for the greatest 
common denominator (equivalent to the algorithm for rewriting a number as a continued fraction).
The idea is to show that there is a string of operators from $\{\mathbf N, \mathbf T\}$  which brings any 
symbol in $\mathcal D$ to the central $0/(1, -1)$. 
Without loss of generality we assume in the following that 
$\dot x$ and $\beta$ are nonnegative.

Since  $\beta\gamma = \dot x^2 -1$,  one of $\beta$ and $\gamma$ must be smaller than $\dot x$.  
Using $\mathbf N$ of (\ref{eq:transfbold}) if necessary, we can make a transition to a new symbol, 
where the curvature entry is the smaller one, say:
$$
\frac{\dot x}{\beta,\gamma}, \qquad \beta<\dot x
$$
Next we can diminish the value of $\dot x$ until it is nonnegative but smaller than $\beta$  
by using $\mathbf T^{-1}$ an appropriate number of times.  A handy formula
for its higher power is:
$$
    	\mathbf T^{-n} : \  \ \frac{\dot x}{\beta, \, \gamma} 
                                 \ \mapsto \  
                                 \frac{\dot x \pm n\beta}{\beta, \  n^2\beta  \pm  2n\dot x+ \gamma}   
$$
Now it is the new co-curvature $\gamma$ that is smaller than $\dot x$. 
We repeat the two steps until the denominator $\dot x$ vanishes.  
This is possible only if the symbol is 
$$
\frac{0}{1,\, -1}
$$
(the ``denominator'' up to sign, adjustable by inversion $\mathbf N$), 
which corresponds to the principal circle at $x = 0$.    

The only obstacle one could encounter while executing this algorithm 
(possibly at the outset) is the $\beta$-entry, the curvature, becoming $0$, 
making   
diminishing of the value of $\dot x$ impossible.  
In this case $\beta\gamma = 0$ enforces, 
by normalization (\ref{eq:norm}), that the ``numerator'' of the circle symbol 
be 1, and $\gamma$ be odd to preserve the (2,1)-parity.  
However, this symbol {\it is} in $\mathcal D$, since it is an image of the following chain of operations:
$$
\frac{0}{1,-1}  
\quad \xrightarrow{\ \ \mathbf T \ \ } \quad
\frac{1}{1,0}  
\quad \xrightarrow{\ \ \mathbf N \ \ } \quad
\frac{1}{0,1}  
\quad \xrightarrow{\ \ \mathbf T^{\,k}\ \ } \quad
\frac{1}{0,2k+1}  
$$
And this ends the proof.  
$\square$
\\


Summarizing, the last theorem may be restated in a more concise but less practical way:

\begin{theorem}[{\rm =}3.7] 
\label{theorem:mainagain}
The circles (including lines) constituting the orbit of the unit circle in the complex plane $\mathbb{C}$ 
via the point-wise M\"obius action of the modular group $\PSL(2,\mathbb Z)$ on $\mathbb{C}$ are 
in one-to-one correspondence with the integer solutions to
$$
                    k^2 - 1 = nm
$$
where $n\geq 0$ and exactly two of the three integers $(m, n, k)$ are odd.  The corresponding circle's symbol is
$$
\frac{k}{n,\, m}.
$$
\end{theorem}

\section{Some patterns}

With an algorithm at hand, we can now draw an arbitrarily large portion of the Dedekind tessellation.  
In particular, we can readily produce the set of all "Dedekind fractions" $k/n$ for a fixed maximum $n$. 
Figure \ref{fig:bigmars} below presents such fractions in graphical form for $k < n < 1000$.

~

\begin{figure}[h]
\centering
\includegraphics[scale=.94]{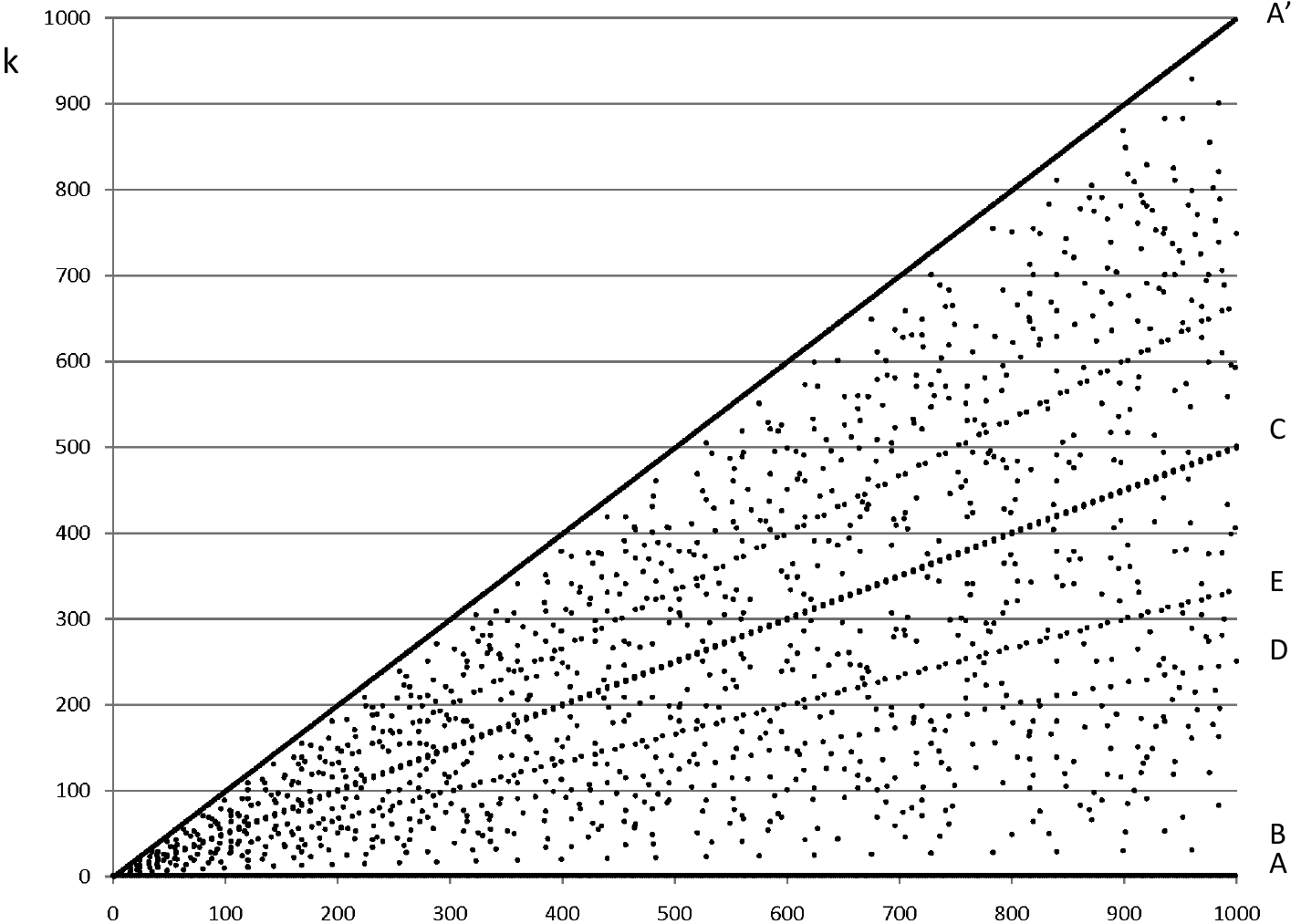}
\caption{\small Dedekind fractions.  A dot at $(n,k)$ represents fraction $k/n$}
\label{fig:bigmars}
\end{figure}

~

A few series are conspicuous. Among them are:
$$
\begin{array}{clrl}
    A\,     &\hbox{Very bottom series: }      &     1/(2\ell+1):     
             &\hbox{\small 1/1, 1/3,  1/5, 1/7, \ldots}\\ 
    A'      &\hbox{Upper diagonal series:}    &   2\ell/(2\ell+1):    
             & \hbox{\small 0/1, 2/3, 4/5, 6/7, \ldots}\\
    B\,     &\hbox{Right above the bottom line: }     &   \ell/(\ell^2-1) :    
             & \hbox{\small 2/3, 3/8, 4/15, 5/24, \ldots}\\
    C\,     &\hbox{In the middle of the cloud:  }    &   (4\ell\pm 1)/8\ell:  
             &\hbox{\small 3/8, 5/8, 7/16, 9/16, \ldots}\\
    D\,     &\hbox{Close to the line of slope 1/4:} &  (8\ell\pm 3)/8(4\ell\pm 1) : 
              &\hbox{\small 5/24, 11/40, 13/56, 19/72, ...}\\
    E\,      &\hbox{Close to the line of slope 1/3:}           &      (6\ell\pm 2)/(18\ell\pm 3) : 
             &\hbox{\small 2/3, 4/15, 8/21, 10/33, \ldots}
\end{array}
$$

Other patterns, some quite visible, emerge as well.

\newpage 

\noindent
{\bf Occurrences of the golden ratio and Fibonacci numbers.}   
Recall that the Fibonacci sequence of numbers $(F_n)$  is defined 
recursively by setting $F_0 = 0$,  $F_1 = 1$,  and   $F_{k+2} = F_k+F_{k+1}$.  
This generates 
$$
\mathbf 0,\,  1,\,  \mathbf 1,\,  2,\, \mathbf  3,\,  5,\, \mathbf 8,\,  13,\, \mathbf 21,\,
                    34,\,  \mathbf 55,\,  89,\,  \mathbf 144, \ldots
$$
(the even-indexed entries are in bold face for quick identification).   
The following fractions are found among the centers of Dedekind circles 
(see Figure \ref{fig:table}): 
$$
\begin{array}{clccl}
A. \ \      
&
\dfrac{1}{1},\;
\dfrac{2}{3},\;
\dfrac{5}{8},\;
\dfrac{13}{21},\;
\dfrac{34}{55},\;
\dfrac{89}{144},\ldots
&
\dfrac{F_{2n-1}}{F_{2n}}\;
&\mapsto
&\dfrac{\sqrt{5}-1}{2} = \tau \equiv \varphi^{-1}  \\
\\
B.  \ \ 
&
\dfrac{1}{1},\;
\dfrac{5}{3},\;
\dfrac{13}{8},\;
\dfrac{34}{21},\;
\dfrac{89}{55},\;
\dfrac{233}{144},\;
&
\dfrac{F_{2n+1}}{F_{2n}}\;
&\mapsto
&\dfrac{\sqrt{5}+1}{2} = \varphi  
\end{array}
$$
Sequence (A) may be obtained by a recurrent application of 
$\mathbf{NT}$ (translation followed by inversion):
$$
A  =  \mathbf N \mathbf T   = 
\left[
\begin{array}{ccc}
1&1&0\\
2&1&1\\
0&1&0
\end{array}\right]
:\quad 
\left[
\begin{array}{ccc}
\dot x\\
\beta\\
\gamma
\end{array}\right]
\mapsto
\left[
\begin{array}{ccc}
\dot x+\beta\\
2\dot x+ \beta+\gamma\\
\beta
\end{array}\right]
$$
The circle symbol consists of three consecutive Fibonacci numbers:
\begin{equation}
\label{eq:FFF}
\dfrac{F_{2n-1}}{F_{2n},\, F_{2n-2}}\; .
\end{equation}
Sequence (B)  results from the same symbol (\ref{eq:FFF}) with the curvature and co-curvature switched. 
Geometrically, it may be obtained from the circles of (A) by inversion in the unit circle or, equivalently,  
by a unit transformation (adding one to the fraction). 
Interestingly, the normalization (\ref{eq:normbig}, \ref{eq:norm}) in the Minkowski space corresponds 
to a well-known identity for Fibonacci numbers:    
$$
(F_{2n-1})^2 - 1  =  F_{2n}F_{2n-2}
$$

The following two sequences of Dedekind circles are also related to Fibonacci numbers:
$$
\begin{array}{clccc}
C.  \ \    
&
\dfrac{1}{3},\;
\dfrac{4}{5},\;
\dfrac{9}{16},\;
\dfrac{25}{39},\;
\dfrac{64}{105}, \ldots
&
\dfrac{F_n^2}{F_{n-1}F_{n+2}}\;
&\mapsto
&\dfrac{\sqrt{5}+1}{2} = \tau  \\
\\
D. \ \
&
\dfrac{4}{3},\;
\dfrac{9}{5},\;
\dfrac{25}{16},\;
\dfrac{64}{39},\;
\dfrac{169}{105},\ldots
&
\dfrac{F_{n}^2}{F_{n-2}F_{n+1}}\;
&\mapsto
&\dfrac{\sqrt{5}-1}{2} = \varphi  
\end{array}
$$
The the circle symbols of $C$ are composed of five consecutive Fibonacci numbers
\begin{equation}
\dfrac{F_{n}^2}{F_{n-1}F_{n+2},\; F_{n+1}F_{n-2}}\;
\end{equation}
Switch the entries in the denominator to get (D).  The related identity corresponding to normalization is
$$
(F_n)^4 - 1 \ = \  F_{n-2} F_{n-1} F_{n+1} F_{n+2}       \, ,
$$
which is a product of two known identities,  $(F_n)^2 \pm 1  =  F_{n-1} F_{n+1}$  
and $(F_n)^2 \mp 1 =  F_{n-2} F_{n+2}$,  
with the sign in front of the $1$ depending on the parity of $n$.

\newpage 
\section{Dedekind circles in the Poincar\'e disk}

The standard map  $\Phi:\ {\mathbf H}\to\mathbf D$ from the Poincar\'e half-plane 
to the Poincar\'e disk will carry any circle in the complex plane with center on the real line 
(thus perpendicular to it) to a circle perpendicular the unit circle, thus with the symbol 
of the form in which the curvature and the co-curvature coincide, 
$\beta=\gamma$ (circles invariant under inversion $\mathbf N$). 
We shall use two notations:
$$
\frac{  \dot x,\,\dot y }{  \beta,\beta} 
\ \equiv \ 
\dotfraction{ \ 2\ddot x,\,2\ddot y \ }{  2\beta, 2\beta} \, .
$$ 
This is to preserve the integrality of the description, for we shall see that some 
of the resulting $\dot x$ $\dot y$ and $\beta$ are half-integers.

\begin{proposition}
The map $\Phi:\ {\mathbf H}\to\mathbf D$ from the Poincar\'e half-plane 
to Poincar\'e disk induces the following  map (and its inverse) of the Dedekind circles: 
$$
 \Phi: \ 
\frac{  \dot x  }{  \beta,\gamma} \ \to \ 
        \dotfraction{ \ 2\dot x,\;  \gamma-\beta \ }{  \beta+\gamma, \beta+\gamma }
\qquad \quad
 \Phi^{-1}: \ 
\dotfraction{  \ddot x, \; \ddot y }{ k, k} \ \to \ 
\frac{  \ddot x /2 }{   k-\ddot y, \; k+\ddot y}
$$
\end{proposition}

~\\
{\bf Proof:}  The map may be represented as a sequence of simple transformations
$$
                  \Phi = V_1   \circ  S_2 \circ N \circ V_{-1}
$$
where $V_p$ denotes vertical translation by vector $p$, 
$S_s$ = dilation by factor $s$, and $N$ = inversion in the central unit circle.  
A simple matrix multiplication, after consulting Proposition \ref{prop:transf}, 
gives the result:
$$
\Phi \ = \  \frac{1}{2}
\left[
\begin{array}{rr|rr}
2 & 0 & 0 & 0 \\
0 & 0 & -1&1 \\ \hline
0&2&1&1 \\
0&-2&1&1
\end{array}  \right]
                                                               \qquad
\Phi^{-1}  \ = \  \frac{1}{2}
\left[ \begin{array}{rr|rr}
2 & 0 & 0 & 0 \\
0 & 0 &1&-1 \\ \hline
0&-2&1&1 \\
0&2&1&1
\end{array} \right]
$$
To eliminate the half-values we rescale the size by half, defining 
$\Phi=S_{1/2} \circ \Phi_0$.
For the Dedekind system of circles $y=0$, which gives the result. 
$\square$
\\

\def\PoincareDisk{
\draw (-6, 0)--(6, 0) ;   
\foreach \a/\b/\c   in {  
2/1/1,   2/3/3,   4/2/4,   2/5/5,   2/7/7,   8/2/8,   2/9/9,  6/7/9,  2/11/11,  10/5/11,  12/2/12,  2/13/13,  2/15/15,
8/14/16,  16/2/16,  2/17/17,  2/19/19,  14/13/19,  20/2/20,  2/21/21,  18/11/21,  2/23/23,  22/7/23,  16/18/24,
24/2/24,  2/25/25,  10/23/25,  2/27/27,  28/2/28,  2/29/29,  22/19/29,  26/13/29,  2/31/31,  26/17/31,  32/2/32,
2/33/33,  2/35/35,  12/34/36,  20/30/36,  36/2/36,  2/37/37,  2/39/39,  30/25/39,  38/9/39,  40/2/40,  2/41/41,
34/23/41,  2/43/43,  22/37/43,  28/34/44,  44/2/44,  2/45/45,  2/47/47,  48/2/48,  2/49/49,  14/47/49,  38/31/49,
46/17/49, 2/51/51,  42/29/51, 2/53/53, 38/37/53, 2/55/55, 50/23/55, 56/2/56, 32/46/56, 2/57/57, 
2/59/59, 26/53/59, 46/37/59, 58/11/59,     60/2/60, 52/30/60, 2/61/61, 58/19/61, 50/35/61, 2/63/63, 42/47/63,
64/2/64, 16/62/64, 40/50/64, 2/65/65, 2/67/67, 68/2/68, 28/62/68, 2/69/69, 54/43/69, 2/71/71, 58/41/71,
72/2/72, 2/73/73, 2/75/75, 70/27/75, 76/2/76, 68/34/76, 44/62/76, 2/79/79, 62/49/79, 80/2/80,
2/81/81, 66/47/81, 18/79/81, 74/33/81,    2/83/83, 82/13/83, 84/2/84, 52/66/84, 2/85/85, 2/87/87, 
88/2/88, 32/82/88, 2/89/89, 86/23/89, 70/55/89, 2/91/91, 74/53/91, 92/2/92, 68/62/92, 2/93/93, 42/83/93,
2/95/95, 96/2/96, 78/56/96, 2/97/97, 98/2/98,     2/99/99, 34/93/99, 54/83/99, 78/61/99, 100/2/100, 20/98/100  }
\draw 
          (\a/ \c, \b/ \c)  circle  (2/\c)
          (\a/ \c, -\b/ \c)  circle  (2/\c)
          (-\a/ \c, \b/ \c)  circle  (2/\c)
          (-\a/ \c, -\b/ \c)  circle  (2/\c) ;
\draw [color=black]  (0, 0) circle (1);
}
%
%
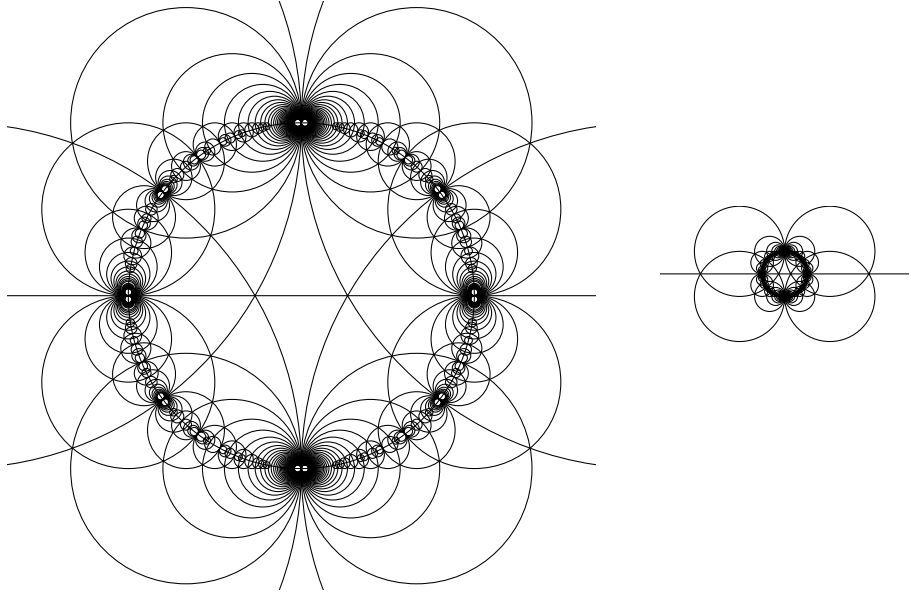
\begin{figure}[h!]
\centering
\begin{tikzpicture}[scale=2.3]
\clip (-1.7,-1.7) rectangle (1.7,1.7);
\PoincareDisk
\end{tikzpicture}
\qquad 
\begin{tikzpicture}[scale=.3]
\clip (-5.5,-14) rectangle (5.5,3);
\PoincareDisk
\end{tikzpicture}

\caption{Dedekind tessellation in the Poincar\'e disk. On the right is a telescopic view}
\label{fig:disk}
\end{figure}

Using this transformation one may easily produce the picture of 
Dedekind circles $\mathcal D$ in the Poincar\'e disk version, $\mathcal D_0$;   
the result is shown in Figure \ref{fig:disk}.  
The symmetry group is now (the projective version of) $\SU_{1,1} \left(\,\mathbb Z[i]\,\right)$, by which we understand 
the modular group conjugated by the matrices that carry the unit disk to the upper-half plane and back:
$$
\PSU_{1,1}(\mathbb Z[i])    
                            =  \left[\begin{matrix}  \ \,  1 &  \!\! \!\!  -i \\ -i&1\end{matrix}\right]
                                   \cdot \PSL(2,\mathbb Z) \cdot
                                 \left[\begin{matrix}1&i \\i&1\end{matrix}\right]
                            \cong \left\{  \left[\begin{matrix}z&\!\! w \\ \overline{w}& \!\!  \overline{z}\end{matrix}\right]
                                  \ \Big| \ z,w\in \mathbb Z[i], \; |z|^2 -|w|^2 = 4 \right\} 
$$
(the conjugating matrices are modified Cayley transforms).
It may be understood as $2\times 2$ integer complex matrices of determinant 4.  
Its projective version $\PSU_{1,1}(\mathbb Z)$ caries $\mathcal D_0$ to itself via M\"obius transformations.
\\

\newpage

The disk version of Theorem \ref{theorem:main} follows: 

\begin{corollary}  
The circles in the orbit of the real axis $\mathbb R\in \mathbb C$ via the induced 
M\"obius action of the modular group $\PSU_{1,1}(\/\/\mathbb Z[i]\/\/)$ on $\mathbb C$ are 
in one-to-one correspondence with integer solutions to
$$
n^2+ 4 = p^2 + q^2
$$
where $n\in \mathbb N$ is either odd or divisible by 4.  In the former case only one of $p$ and $q$ is even, say $p$. 
In the latter case exactly one of them is divisible by 4 (although both are even), say $p$.  
In both cases the circle's symbol is 
$$
          \dotfraction{ \ p,q \ }{ n, n}   
                                \qquad \rightarrow 
                                         \quad \hbox{\sf circle at } \ (p/n, q/n) \ \hbox{\sf with radius} \ R=2/n  
$$ 
The correspondence is not unique only in the case of $n = 0$, where both solutions, 
$(p,q)=(0,\pm 2)$,  represent the same line, $y=0$.
\end{corollary}

The first few symbols ordered by increasing curvature are:
\begin{figure}[H]
\centering
\includegraphics[scale=.9]{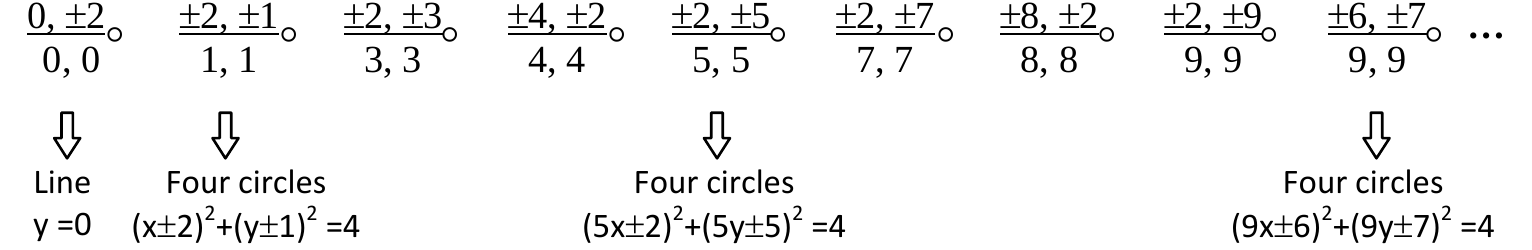}  
\end{figure}
\noindent
Note that at least one circle for each permissible $n$ exists (with 2 for one of the numbers in the ``numerator'').
\\

As before, the system is self-symmetric:  inversion in any of the circles of $\mathcal D_0$ maps the whole system to itself.
The image of the upper half plane constitutes the exterior of the unit circle; 
its interior is the image of the lower-half plane.

Interestingly, the pattern 
has additional meaning:  It consists of circles that are among the inversive symmetries of the Apollonian window,
a special case of the Apollonian disk packing \cite{jk-Descartes}.  
Similarly, the Dedekind tessellation of the upper half plane contains symmetry circles 
for the unbounded Apollonian disk packing (the Apollonian belt, or Ford circles).

\section{Code}

The results presented may be of some value for computer graphics.  
The two figures below present algorithms
that allow one to draw both Dedekind tessellations 
in a consistent way, following a monotonic order of curvatures.
\\


\begin{figure}[h!]
{\sf\small
\noindent
\hspace*{.3in}\rule{4.7in}{1.2pt}       
\\[-2pt] \hspace*{.4in}          Choose $N$  (the greatest curvature desired) 
\\[-2pt] \hspace*{.4in}         LINE $x = 1/2$                                                       
\\[-2pt] \hspace*{.4in}          For $n =1$ to $N$                   
\\[-2pt] \hspace*{.6in}               For $k=0$ to $n$           
\\[-2pt] \hspace*{.8in}                    $m=(k^2 - 1)/n$      
\\[-2pt] \hspace*{.8in}                                  If  (\,$n$ = odd  and $m$ = integer\,) \ or \ (\,$8\,|\,n$ and $m$ = odd\,)   
\\[-2pt] \hspace*{1.4in}                                            then CIRCLE at $(k/n, 0) \in \mathbb R^2$ with radius $R = 1/n$  
\\[-2pt] \hspace*{.6in}                Next $k$    
\\[-2pt] \hspace*{.4in}            Next $n$             
\\[-2pt] \hspace*{.4in}            Redraw the above under any desired number of unit horizontal translation  
\\[-8pt] \hspace*{.3in}\rule{4.7in}{1.2pt}
}
\caption{Drawing the Dedekind tessellation in the Poincar\'e half-plane}
\label{fig:alg1}
\end{figure}


\begin{figure}[h!]
{\sf\small
\noindent
\hspace*{.3in}\rule{4.7in}{1.2pt}
\\[-2pt] \hspace*{.4in}          LINE $y = 0$
\\[-2pt] \hspace*{.4in}          For $n =1$ to infinity (or less)
\\[-2pt] \hspace*{.6in}           (If $n/2$ is an odd integer then go to next $n$)
\\[-2pt] \hspace*{.6in}           For $p = 1$ to $n$
\\[0pt] \hspace*{.8in}                 $q = \sqrt{n^2 - p^2 +4}$  
\\[0pt] \hspace*{.8in}                CASES
\\[-2pt] \hspace*{.8in}                 CASE 1:  If $n$ is odd and $q\in\mathbb Z$ and $2\,|\,q$ 
\\[-2pt] \hspace*{1.4in}                              then CIRCLE at $(q/n,\, p/n) \in \mathbb R^2$ with radius $R = 2/n$
\\[-2pt] \hspace*{.8in}                CASE 2:  If $4\,|\,n$ and  $q \in\mathbb Z$ and $4\,|\,q$ 
\\[-2pt] \hspace*{1.4in}                              then CIRCLES at $(\pm p/n, \pm q/n) \in \mathbb R^2$, with radius $R = 2/n$
\\[-2pt] \hspace*{.6in}            Next $p$
\\[-2pt] \hspace*{.4in}            Next $n$
\\[-8pt] \hspace*{.3in}\rule{4.7in}{1.2pt}
}
\caption{Drawing the Dedekind tessellation in the Poincar\'e disk}
\label{fig:alg2}
\end{figure}

~\\\\
{\bf Acknowledgments:}  
The author thanks Philip Feinsilver for invaluable remarks and comments.



\begin{thebibliography}{9}

\bibitem 
{Ded} Richard Dedekind, 
        Schreiben an Herrn Borchard \"uber die Theorie der elliptische Modulfunktionen. 
        {\it J. Reine Angew. Math.}  Bd. {\bf 83}  (1877), 265-292.


\bibitem
{Kn} Felix Klein, 
        Uber die Transformation der elliptischen Funktionen und die Aufl\"osung der Gleichungen f\"unften Grades. 
       {\it Math. Ann.}   Bd. {\bf 14} (1879), 111-172. 

\bibitem
{jk-Descartes} Jerzy Kocik, 
        A theorem on circle configurations. arXiv:0706.0372. 

\bibitem
{jk2} Jerzy Kocik, 
        Clifford algebras and Euclid's parameterization of Pythagorean triples. 
          {\it Adv. Appl. Clifford Al.} {\bf 17} (2007),  71-93.


\bibitem
{LMW} 
Jeffrey C. Lagarias, Colin L. Mallows, Allan R. Wilks,
        Beyond the Descartes circle theorem. {\it Am. Math. Monthly} 109 (2002), 338--361.

\bibitem
{LeB} Lieven LeBruyn, Dedekind or Klein?  \hfil\break
          \url{www.neverendingbooks.org/dedekind-or-klein}

\bibitem
{Sti} John Stillwell, 
        Modular Miracles. {\it Am. Math. Monthly}  {\bf 108} (2001), 70-76.


\end{thebibliography}
\end{document}